\documentclass[12pt]{article}
\usepackage{latexsym,amsfonts}
\usepackage{picinpar}
\newcommand{\fpr}{\hfill $\Box$\\}
\newtheorem{theorem}{Theorem}
\newtheorem{proposition}{Proposition}

\font\tenmsb=msbm10 scaled\magstep1\font\sevenmsb=msbm7 scaled\magstep1
\font\fivemsb=msbm5 scaled\magstep1\newfam\msbfam
\textfont\msbfam=\tenmsb\scriptfont\msbfam=\sevenmsb
\scriptscriptfont\msbfam=\fivemsb
\def\G{\Gamma}
\def\la{\lambda}
\def\E{\textsf{E}}
\def\P{\textsf{P}}
\def\Pr{\textsf{Pr}}
\def\var{\textsf{var}}
\def\Q{\textsf{Q}}
\newenvironment{abstract1}{%
         \small
     \begin{center}  
{\bfseries {Abstract}}%
               \end{center}}

\begin{document}

\title{Random Partitions with non negative $r^{th}$ differences}
\author{Rod Canfield$^{a}$, Sylvie Corteel$^{b}$, Pawel Hitczenko$^{c}$\\
\vspace{.5cm}\texttt{erc@cs.uga.edu}, \texttt{syl@prism.uvsq.fr}, 
\texttt{phitczen@mcs.drexel.edu}\\
$^{a}$ Department of Computer
Science, \\
University of Georgia Athens, GA 30602, USA\\
$^{b}$  PRiSM, Universit\'e de Versailles,\\
45 Av. des Etats Unis, 78035 VERSAILLES, France\\
$^c$  Department of Maths and Computer
Science, \\
Drexel University, Philadelphia, PA 19104, USA}

\date{\today}
\maketitle

\begin{abstract1}
Let $P_r(n)$ be the set of partitions of $n$ with non negative $r^{th}$ 
differences.
Let $\lambda$ be a partition 
of an integer $n$ chosen uniformly at random among the set $P_r(n)$
Let $d(\lambda)$ be a positive $r^{th}$ difference
chosen uniformly at random  in $\lambda$. The aim of this work 
is to show that for every $m\geq 1$, the probability that
$d(\lambda)\geq m$ approaches the constant $m^{-1/r}$ as $n\rightarrow\infty$
This work is a generalization of a result on integer partitions 
\cite{cps} and was motivated by a recent identity by Andrews, 
Paule and Riese's Omega package \cite{apr98}.
To prove this result we use bijective, asymptotic/analytic
and probabilistic combinatorics.
\end{abstract1}

\vspace{.25in}

\section{Introduction}

Let us first start by a few definitions and notations.
A partition $\lambda$ of $n$ 
is a sequence of integers

\[
\lambda=(\lambda_1,\lambda_2,\ldots,\lambda_k)\ \ {\rm with}\  
\lambda_1\geq \lambda_2\ldots \ge \lambda_k\ge 1\ \ {\rm and} \ 
\sum_{j=1}^k \lambda_j=n.
\]
The same partition $\lambda$ can also be written in its frequential
notation, that is~: 
$$\lambda=n^{m_n}(n-1)^{m_{n-1}}\ldots1^{m_1}
\ \ {\rm with}\ m_j=|\{i\ |\ \lambda_i=j\}|,\  1\leq j\leq n.$$ 
The number $m_j$ is called the multiplicity of the part 
$j$ in $\lambda$. In the sequel we will use both representations.
We now define the $r^{th}$ differences.
Let  $\lambda=(\lambda_1,\lambda_2,\ldots,\lambda_k)$ be a partition 
of $n$ and $\Delta^r(\lambda)=(\Delta^r_1(\lambda),\ldots,\Delta^r_k(\lambda))$ be its $r^{th}$
differences. The $r^{th}$ differences can be computed by the following recurrence~:
\[
\Delta^r_i(\lambda)=\left\{
\begin{array}{ll}
\lambda_{i}&{\rm if}\ i=k\  {\rm or}\ r=0\\
\Delta^{r-1}_i(\lambda)-\Delta^{r-1}_{i+1}(\lambda)&{\rm otherwise}
\end{array}\right. 
\]
In the sequel we will write $\Delta^r_i(\lambda)$~: $\Delta^r_i$ for short.
Let $P_r$ be the set of partitions with non negative 
$r^{th}$ differences, that is to say, 
$(\lambda_1,\lambda_2,\ldots,\lambda_k)\in P_r$
if and only if $\Delta^r_i\geq 0$ for $1<i< k$.
Let $P_r(n)$ be the set of partitions $\lambda\in P_r$ with $|\lambda|=n$
and let $p_r(n)=|P_r(n)|$.\\

This work was motivated by two previous results. 
The first result is an identity on partitions with non negative
$r^{th}$ differences. It was discovered by Andrews, 
Paule and Riese's Omega package.

\begin{theorem}{\rm \cite{an98,apr98}}
There is a one-to-one correspondence between the partitions in $P_r(n)$ and
the partitions of $n$ into parts
${i+r\choose r}$, $i\ge 0$.
\label{omega}
\end{theorem}

The second result is on
ordinary integer partitions (partitions with non negative $1^{st}$ 
differences)~:
\begin{theorem}{\rm \cite{cps}}
Let $\lambda$ be a partition
of an integer $n$ chosen uniformly at random among the set of all
partitions of $n$.
Let $d(\lambda)$ be a part size
chosen uniformly at random from the set of all part sizes
that occur in $\lambda$. For every $m\geq 1$, the probability that
$d(\lambda)\geq m$ approaches the constant $1/m$ as $n\rightarrow\infty$.
\label{ordi}
\end{theorem}

Note that of $r=1$ there exists a one to one correspondence between 
the number of positive $1^{st}$ differences in any given partition
and the number part sizes of its conjugate. Our aim is therefore
to generalize the result of Theorem \ref{ordi} by using the identity of 
Theorem \ref{omega}.\\


Let us now state our generalization~:

\begin{theorem}
Let $\lambda$ be a partition
of an integer $n$ chosen uniformly at random among the set $P_r(n)$ of all
partitions of $n$ with  non negative $r^{th}$ differences.
Let $d(\lambda)$ be a positive $r^{th}$ difference
chosen uniformly at random from the set of all positive differences
that occur in $\lambda$. For every $m\geq 1$, the probability that
$d(\lambda)\geq m$ approaches the constant $m^{-1/r}$ as $n\rightarrow\infty$.
\label{main}
\end{theorem}

The purpose of this paper is to prove Theorem \ref{main}.
We now present the organization of the paper.
We will first give in Section 2 
a simple bijection of the identity of Theorem \ref{omega} that 
gives several refinements of the identity. 
This bijection was 
advertised/announced by Zeilberger in his very own journal \cite{co98}.
We then present in Section 3 some asymptotic results
on the number of partitions in $P_r(n)$ and, thanks to the bijection,
on the average number 
of positive $r^{th}$ differences in the partitions in $P_r(n)$.
We finally use in Section 4
some probabilistic arguments which generalize the works
on ordinary partitions \cite{f,cps}. The association of the three parts 
gives us the proof of our result. We conclude the paper by presenting 
some future work in Section 5.

\section{Bijective Combinatorics}

In this section we are going to present a bijection $f$ between
the partitions in $P_r(n)$ and the partitions of $n$ into parts 
${i+r\choose r}$, $i\ge 0$.
Let
$\lambda=(\lambda_1,\lambda_2,\ldots,\lambda_k)$ be a partition in  $P_r(n)$
then its image by the bijection $f$ in its frequential notation is~:
\[
f(\lambda)={{k-1+r}\choose{r}}^{\Delta^r_k}{{k-2+r}\choose{r}}^{\Delta^r_{k-1}}\ldots{r\choose r}^{\Delta^r_1}
\]
where the ${\Delta^r_i}$ ($1\le i\le k$) are the $r^{th}$ differences
of $\lambda$. Is is clear that $f(\lambda)$  
is a partition into parts
${r+i\choose r}$ with $i\geq 0$.
Now let us prove that $f(\lambda)$  
is a partition of $n$.
>From the definition of $\Delta^r(\lambda)$,
it is easy to see that
\[
\lambda_i=\sum_{j=i}^k {r+k-j-1 \choose r-1} \Delta^r_{j}.
\]
As ${r+i-1\choose r}=\sum_{j=0}^{i-1}{r+j-1\choose r-1}$, 
we have $|f(\lambda)|=|\lambda|$. 
We can reconstruct $\lambda$ from $\Delta^r(\lambda)$. The reverse mapping
$f^{-1}$ is then easy to define. Let $\mu$ be a partition
of $n$ into parts ${r+i\choose r}$ with $\mu_1={k+r-1\choose r}$. 
Let $\mu^{(i)}$ be the multiplicity
of the part ${i-1+r\choose r}$ in $\mu$, $1\leq i\leq k$. 
Then
\[
f^{-1}(\mu)=\left(\sum_{j=1}^k {r+k-j-1 \choose r-1}\mu^{(j)},
\sum_{j=2}^k {r+k-j-1 \choose r-1}\mu^{(j)},\ldots
,{r-1 \choose r-1}\mu^{(k)}\right).
\]
It is then easy to see that $f$ is a bijection.
Note that if $r=1$ then $f(\lambda)=\lambda'$, the conjugate of $\lambda$.\\

As we said in the introduction, this bijection gives some refinements 
of the identity.
Let us now present these refinements~:
\begin{theorem}
There is a one-to-one correspondence between the partitions 
in $P_r(n)$ with $k$ parts and $j$ positive
$r^{th}$ differences and
the partitions of $n$ into parts ${i+r\choose r}$, $i\ge 0$, whose largest part is ${k-1+r\choose r}$ and with
$j$ parts sizes.
\end{theorem}
\noindent{\bf Proof.} Straightforward with the bijection.\fpr

Let us now illustrate the refinements thanks to generating functions
and a recurrence.
Let $p_r(n,k,m)$ be the number of partitions in $P_r(n)$ 
with $k$ parts and $\sum_{i=1}^{k} \Delta^r_{i}=m$.
Then 
\[
\sum_{n,k\geq 0} p_r(n,k,m)q^ny^kx^m=1+\sum_{k\geq 1}
\frac{y^k x q^{{k-1+r}\choose{r}}}
{(1-xq^{r\choose r})(1-xq^{r+1\choose r})\ldots
(1-xq^{{r+k-1}\choose{r}})}
\]

Let $p_{r,\leq}(n,k,m)$ be the number of partitions in $P_r(n)$ 
with at most $k$ parts and $\sum_{i=1}^{k} \Delta^r_{i}=m$, then
\[
\sum_{n\geq 0} p_{r,\leq}(n,k,m)q^nx^m
=\prod_{i=0}^k (1-xq^{{i+r}\choose{r}})^{-1}.
\]

Let $p_{r,\leq}(n,k)$ be the number of partitions in $P_r(n)$ with at most $k$ parts. 
We can get an easy recurrence to compute this number~:

\[
p_{r,\leq}(n,k)=\left\{ \begin{array}{ll}
0&{\rm if}\ n<0\ {\rm or}\ k=0\ {\rm and}\ n>0\\
1&{\rm if}\ n=0\ {\rm and}\ k=0\\
p_{r,\leq}(n-{{r+k-1}\choose{r}},k)+p_{r,\leq}(n,k-1)&{\rm otherwise}.
\end{array} \right.
\]

\section{Asymptotic Combinatorics}

\def\B{E}
\def\oh{{\rm o}}

In this section we develop an asymptotic formula for
$p_S(n)$, the number of partitions of the integer
$n$ whose parts lie in a set $S$ which is the image
of a given polynomial.   As a special case, we have
for our $p_r(n)$:

\begin{theorem}
As  $n \rightarrow \infty$, we have:
\begin{equation}
p_r(n)\sim c n^{-(r+1)/2} \exp\{ C(1+r)n^{1\over 1+r} \}
\end{equation}
where the constants $C$ and $c$ are given by
\begin{equation}
\label{Ceq}
C = \bigl\{ r!^{1/r} r^{-1} \zeta(1+r^{-1}) \Gamma(1+r^{-1})
\bigr\}^{r \over r+1}
\end{equation}
and
\begin{equation}
c = \bigl( {C \over 2\pi} \bigr)^{r+1 \over 2}
(r!)^{-r/2} (1+r^{-1})^{-1/2} \, \prod_{j=0}^{r-1} j!.
\end{equation}
\label{asymp1}
\end{theorem}

To prove Theorem \ref{asymp1} we can
take $P(x)={x+r-1\choose r}$ in
Theorem \ref{asymp2} below.  It is seen that
$d=r$, $A=1/r!$, $B=1/(2(r-2)!)$, and
$\rho_j = j-1$.  We now
state and prove the more general result.

\begin{theorem}
Let $P(x) = Ax^d + Bx^{d-1}+\cdots$ be a
polynomial of degree $d$ which is positive and
increasing for $x \ge 1$,
and which assumes integral values for integral $x$.  Let
$S$ be the set $\{P(1), P(2), \dots \}$, and assume that
the set $S$ has gcd $1$.
Let $\rho_j$ be the negatives of
the roots of $P(x)$, so that
\[
P(x) = A \prod_{j=1}^d (x+\rho_j)
\]
(in particular, $B/A=\sum\rho_j$). Then,
$$ p_S(n) ~~=~~ c \, n^{-\kappa-{1 \over 2}} 
\, \exp\{ (1+d)Cn^{1 \over d+1} + o(1) \},$$
where
\begin{eqnarray*}
\kappa &=& {Ad+B \over A(d+1)}  \\
C      &=& \bigl\{ A^{-1/d} d^{-1} \zeta(1+d^{-1})
\Gamma(1+d^{-1}) \bigr\}^{d \over d+1}  \\
\end{eqnarray*}
and
\begin{equation}
c = \prod_{j=1}^d \Gamma(1+\rho_j) \, C^{1+B/Ad} \,
A^{{1 \over 2}+B/Ad} \, (1+d^{-1})^{-1/2}
\, (2\pi)^{-(d+1)/2}
\end{equation}
\label{asymp2}
\end{theorem}

\noindent {\bf Proof.} Let $N(u)$ be the counting function
associated with the set $S$:
$$ N(u) ~~=~~ \#\{h \in S: h \le u\}.$$
Under the assumptions that
$$ N(u) ~~=~~ \B u^{\beta} ~+~ R(u) $$
with $\B, \beta > 0$ and
$$ \int_0^u {R(v) \over v}dv ~~=~~ b_1 \log u ~+~ b_2 ~+~ \oh(1),
~~~ u \rightarrow \infty $$
Ingham \cite{in41} has shown that for
$n \rightarrow \infty$
$$                                            
{\hat{P}(n) - \hat{P}(n-h)  \over h}
~~\sim~~
\bigl({1-\alpha \over 2 \pi } \bigr)^{1/2}  \,
e^{b_2} \, M^{-(b_1-1/2)\alpha} \, n^{(b_1-1/2)(1-\alpha)-1/2}
\, e^{\alpha^{-1}(Mn)^{\alpha}},
$$
for each fixed $h \in S$, where $\alpha = \beta/(1+\beta)$,
$$ M ~~=~~ \bigl\{\B\beta\Gamma(1+\beta)\zeta(1+\beta)\bigr\}^{1/\beta},$$
and
$$ \hat{P}(n) ~~=~~ \sum_{j \le n} p_S(j). $$

Observe that the right side of Ingham's formula is independent of
$h$; we shall denote it $G(n)$. Observe further that for any fixed
$k$ we have $G(n-k)\sim G(n)$; hence, for any fixed
integer $k$ and fixed $h \in S$,
\begin{equation}
p(n-k)+p(n-k-1)+\cdots p(n-k-h+1) ~~\sim~~ h G(n).
\label{IngSum}
\end{equation}
Because the gcd of the set $S$ is by assumption 1, there must
be two finite, disjoint subsets $S_1, S_2 \subseteq S$, and
for each $h \in S_1 \cup S_2$ a positive integer $r_h$ such
that
$$ \sum_{h \in S_1}r_h h ~~=~~ 1 ~+~ \sum_{h \in S_2}r_h h. $$
Let $H=\sum_{h \in S_2}r_h h$.  If we divide the integers in
the interval 
$[n-H,n-1]$ (there are $H$ such integers) into 
$\sum_{h \in S_2}r_h $ disjoint subintervals, there being
$r_h$ subintervals of size $h$, for $h \in S_2$, and we
apply (\ref{IngSum}) to each of these, we obtain
$$ \sum_{j=1}^H p_S(n-j) ~~\sim~~ \bigl( \sum_{h\in S_2} r_h h \bigr)
\, G(n), ~~ n \rightarrow \infty. $$
Similarly, we may divide the integers in the interval $[n-H,n]$
(there are $H+1$ such integers) into
$\sum_{h \in S_1}r_h $ disjoint subintervals, there being
$r_h$ subintervals of size $h$, for $h \in S_1$, and find
$$ \sum_{j=0}^H p_S(n-j) ~~\sim~~ \bigl( \sum_{h\in S_1} r_h h \bigr)
\, G(n), ~~ n \rightarrow \infty. $$
Subtracting,
$$ p_S(n) ~~\sim~~ G(n).  $$
It may have been an oversight on Ingham's part not to
notice that his formula gave an asymptotic for $p_S(n)$ in
the typical situation that $S$ has gcd $1$.  This was also overlooked
by some later writers who reference Ingham's Tauberian
theorem, but noticed, and explained somewhat differently than we have,
in \cite{HH}.

\vskip 3pt

It remains only for us to determine the constants $\B$, $\beta$,
$b_1$, and $b_2$ for the case that $S=P(\{1,2,\dots\})$.  Our
theorem then follows by using Ingham's $G(n)$.  Given that the polynomial
$P$ is positive and increasing, we have
$$ N(P(x)) ~~=~~ [x], $$
the brackets on the right denoting greatest integer.  Hence,
letting $x_*$ denote the unique real solution of $P(x_*)=u$,
\begin{eqnarray*}
\int_0^u{N(v)\over v}dv &=& \int_{\lambda_1}^u {N(v)\over v}dv \\
&=& 
\int_1^{x_*} {[x] P'(x) \over P(x)}dx \\
&=& 
\int_1^{x_*} {x P'(x) \over P(x)}dx 
~-~ {1 \over 2} \int_1^{x_*} { P'(x) \over P(x)}dx
~-~  \int_1^{x_*} {(x-[x]-1/2) P'(x) \over P(x)}dx  \\
&=&
I_1 ~-~ {1 \over 2} \log P(x_*) ~+~ {1 \over 2} \log P(1) ~-~ I_2. \\
\end{eqnarray*}
For the integral $I_1$ we integrate $x/(x+\rho_j)$, and sum on $j$
to find
$$ I_1 ~~=~~ d(x_*-1) ~-~ \sum_{j=1}^d \rho_j \log(x_* +\rho_j)
~+~ \sum_{j=1}^d \rho_j \log(1+\rho_j). $$
Using the estimates, valid for $u \rightarrow \infty$,
$$ d x_* ~~=~~ d (u/A)^{1/d} ~-~ \sum_{j=1}^d \rho_j ~+~ \oh(1) $$
$$ \log(x_*+\rho_j) ~~=~~ \log(x_*) ~+~ \oh(1)
~~=~~ {1 \over d}\log(u/A) ~+~ \oh(1), $$
we then find
$$
I_1 ~~=~~ d (u/A)^{1/d} ~-~  \sum_{j=1}^d \rho_j 
~-~ d ~-~ {1 \over d} \log(u/A) \sum_{j=1}^d \rho_j 
~+~ \sum_{j=1}^d \rho_j \log(1+\rho_j)  ~+~ \oh(1). 
$$

According to the Euler-Maclaurin summation formula,
$$
- \, \int_1^N {x-[x]-1/2 \over x+\rho} dx
~~=~~  \int_1^N \log( x+\rho) dx ~-~ \sum_{k=1}^N\log(k+\rho)
~+~ {1 \over 2} (\log(1+\rho) + \log(N+\rho)) 
$$
We can sum $\log(k+\rho)$ over $k$ with
the Gamma function,
$$ \sum_{k=1}^N\log(k+\rho) ~~=~~ \log\Gamma(N+1+\rho) ~-~ 
\log\Gamma(1+\rho),  $$
and Stirling's formula,
$$ \log\Gamma(N+1+\rho) ~~=~~ (N+\rho)\log(N) ~-~ N ~+~ {1 \over 2} 
\log(2\pi N) ~+~  \oh(1).  $$
Using this and
the definite integral
$$ \int_1^N \log(x+\rho) dx ~~=~~ (N+\rho) \log(N+\rho) ~-~ N ~+~ 1 
~-~ (1+\rho) \log(1+\rho), $$
we conclude, for $N \rightarrow \infty$, 
$$- \, \int_1^N {x-[x]-1/2 \over x+\rho} dx
~~=~~ 1 ~+~ \rho ~-~ {1 \over 2} \log(2\pi) ~+~ \log\Gamma(1+\rho)
~-~ ({1 \over 2}+\rho)\log(1+\rho) ~+~ \oh(1). $$
Summing over the negative roots $\rho_j$ of $P(x)$, we have
$$ -I_2 ~~=~~ d ~+~ \sum \rho_j ~-~ {d \over 2}\log(2\pi)
~+~ \sum\log\Gamma(1+\rho_j) ~-~ {1 \over 2} \log {P(1) \over A}
~-~ \sum \rho_j \log(1+\rho_j) ~+~ \oh(1). $$
Hence, replacing $P(x_*)$ by $u$ and substituting for
$I_1-I_2$, 
$$ \int_0^u{N(v)\over v}dv ~~=~~ d(u/A)^{1/d}
~-~ {B \over Ad}\log(u/A) ~-~ {1 \over 2}\log(u)
~-~ {d \over 2}\log(2\pi) ~+~ \sum\log\Gamma(1+\rho_j) ~+~ \oh(1). $$
From this we can read off, directly,
\begin{eqnarray*}
\B  &=& A^{-1/d} \\
\beta  &=& 1/d \\
b_1 &=& -{1 \over 2} ~-~ {B \over Ad} \\
b_2 &=& \bigl({1 \over 2}+{B \over Ad}\bigr) \log(A) ~-~ {d \over 2}\log(2\pi)
~+~ \sum\log\Gamma(1+\rho_j). \\
\end{eqnarray*}
The proof is then completed, as explained earlier, by substituting these
four values into Ingham's formula.   \fpr

Next, we use this result on $p_r(n)$ to
compute the asymptotic behavior of the average number
of positive $r^{th}$ differences of the partitions in $P_r(n)$. 
This number $\delta_{r}(n)$ can be defined as follows for any $n$~:
$
\delta_r(n)=
\frac{1}{p_r(n)}\sum_{\lambda\in P_r(n)}|\{i\ |\ \Delta^r_i(\lambda)
>0\}|
$. 
Thanks to the bijection it is straightforward to compute this value for any $n$
and indeed~:
\[
\delta_r(n)=\frac{1}{p_r(n)}\sum_{i\ge 0}p_r\left(n-{{r+i}\choose{r}}\right).
\]

\begin{proposition}
For suitable 
constant $A$
\begin{equation}
\label{mean}
\delta_{r}(n)\sim  A n^{1\over 1+r}
\end{equation}
namely, $A$ is given by
\[
A = \Gamma(1+r^{-1}) \, r!^{1/r} \, C^{-1/r}
\]
where the constant $C$ is given in (\ref{Ceq}).
\label{aver}
\end{proposition}

\noindent{\bf Proof. } Since
$$ (n-K)^{{1 \over r+1}} ~~=~~
n^{{1 \over r+1}} ~-~ {1 \over r+1} K n^{{1 \over r+1}-1}
~+~ \oh(1), $$
uniformly for $K = \oh(n^{1-{1\over 2(r+1)}})$,
we may use our asymptotic formula for
$p_r(n)$ to conclude
$$
{p_r(n-K) \over p_r(n)} ~~=~~ \exp\{ -CKn^{{1 \over r+1}-1}\}, $$
as $n \rightarrow \infty$, uniformly for $K$ restricted
as above.  For $K$ we shall take ${r+j-1 \choose r}$,
$j \ge 1$; those terms for which $K$ is out of range
contribute negligibly to the sum; there results:
$$ \delta_r(n) ~~=~~ \sum_{j \ge 1} \exp\{
-C_n {r+j-1 \choose r}
                                 +\oh(1)\}, $$
where $C_n = C n^{{1 \over r+1}-1}$.
We use Mellin's formula
$$ e^{-y} ~~=~~ {1 \over 2\pi i} \, \int_{a-i\infty}^{a+i\infty}
y^{-s} \Gamma(s) ds, $$
which is valid for Re$(y)>0$ and $a>0$, with $y$
replaced by $C_n {r+j-1 \choose r}$.
This gives
$$ \delta_r(n) ~~\sim~~ {1 \over 2\pi i} \, \int_{a-i\infty}^{a+i\infty}
C_n^{-s} D(s) \Gamma(s) ds, $$
where $D(s)$ is the Dirichlet series
$$ D(s) ~~=~~ \sum_{j \ge 1} {r+j-1 \choose r}^{-s}. $$
By using
$$ \{(r+j-1)(r+j-2)\cdots j\}^{-s} ~-~ j^{-rs}  ~~=~~
(-s) \, \int_{j^r}^{(r+j-1)_r} u^{-s-1}du, $$
and some crude estimates, we can prove that $D(s)$ is equal to
$(r!)^s\zeta(rs)$ plus a function analytic at least in the
domain Re$(s)>0$.   Thus, $ \delta_r(n)$ equals
asymptotically the residue of $ C_n^{-s} (r!)^s \zeta(rs) \Gamma(s)$
at $s=1/r$.  Since
$$ \zeta(rs) ~~=~~ {1 \over rs-1} ~+~ \cdots, $$
the latter residue is $C_n^{-1/r} (r!)^{1/r} r^{-1}\Gamma(r^{-1})$,
and we arrive at the stated formula.  \fpr

Let us now compute the asymptotic behavior of the average number of  $r^{th}$ differences greater or equal to $m$ in the partitions in $P_r(n)$. 
This number $\delta_{r,m}(n)$ can be defined as follows for any $n$~:
$\frac{1}{p_r(n)}\sum_{\lambda\in P_r(n)}|\{i\ |\ \Delta^r_i(\lambda)
\geq m\}|$.
Thanks to the bijection it is straightforward to compute this value for any 
$n$ and $m\ge 1$~:
\[
\delta_{r,m}(n)=\frac{1}{p_r(n)}\sum_{i>0}p_r\left(n-m{{r+i}\choose{r}}\right).
\]

\begin{proposition}
\[
\delta_{r,m}(n)\sim m^{-1/r} A n^{1\over 1+r}
\]
\label{averm}
\end{proposition}
\noindent{\bf Proof.} Similar to the proof of Proposition 1. \fpr

A straightforward consequence of Propositions \ref{aver} and \ref{averm}
is then~:
\begin{proposition}
\[
\lim_{n\rightarrow \infty}\frac{\delta_{r,m}(n)}
{\delta_{r}(n)}=\frac{1}{m^{1/r}}
\]
\label{frac}
\end{proposition}
 
\section{Probabilistic combinatorics}

In this section we will complete the proof of Theorem~\ref{main}. In order to
avoid unnecessary distractions, throughout this section we let $K$  denote generic
constant whose 
 value may change from one use
to the next. These various values may depend on $r$, but are completely
irrelevant from the point of view of the asymptotic results. 

The probability in question is the average value of the ratio
\begin{equation}\E\frac{D_{n,r,m}}{D_{n,r}},
\label{ratio}\end{equation}
where $D_{n,r,m}=\sum_kI(\Delta_k^r\ge m)$ and $D_{n,r}=D_{n,r,1}$.
To compute this average we find it convenient to consider the image of $P_r(n)$ under the bijection $f$. Its explicit form tells us that $\Delta_k^r$'s are multiplicities of parts in the partitions of $n$ into parts whose part sizes are in the set $S=\{{i+r\choose r},\,  i=0,1\dots\}$. Thus  (\ref{ratio}) asks for the asymptotic value of the probability that a randomly chosen part size in a random partition of $n$ into parts in the set $S$ has multiplicity at least $m$.
To compute this probability
we will use arguments which are
the generalization of the work on ordinary partitions \cite{cps,f}.
For convenience we will identify $P_r(n)$ with its image under $f$.
Let $\Q$ be the uniform probability measure on the set $P(n)$ of \textrm{all} partitions of $n$ and $\Pr$ the uniform probability measure on $P_r(n)$.  Since $P_r(n)\subset P(n)$ the measure $\Pr$ is a restriction of $\Q$ to $P_r(n)$ or, in other words it is the conditional measure on {\sl all} partitions given that a partition is in $P_r(n)$.
That is, for any subset $A$ of $P_r(n)$
\[\Pr(\la_{(r)}\in A)=\Q(\la\in A\big|\la\in P_r(n))=\frac{\Q(\la\in A\cap\la\in P_r(n))}{\Q(\la\in P_r(n))},\]
where $\la_{(r)}$ signifies that a partition $\la$ is considered as an element of $P_r(n)$.
We are going to use the result of Fristedt \cite{f},  who proved that if
$\{\G_k:\ k\ge1\}$ are independent geometric random variables with the parameters $1-q^k$, respectively, defined on a probability space with the measure $\P$, then, regardless of the value of $q$,  the joint distribution of multiplicities of parts $(m_1,m_2,\dots)$ in a randomly chosen partition of $n$ is equal to that of $(\G_1,\G_2,\dots)$ conditioned on the event $\{\sum_jj\G_j=n\}$.
Hence,
\begin{eqnarray*}
\Q(\la\in P_r(n))&=&\P((\G_j)\in P_r(n)\Big|\sum j\G_j=n) \\
&=&\frac{\P(\{\G_j=0,\ j\ne{\ell+r\choose r}\}\cap\{\sum j\G_j=n\})}
{\P(\sum j\G_j=n)}
 \\
&=&\frac{\P(\{\G_j=0,\ j\ne{\ell+r\choose r}\}\cap
\{\sum_\ell{\ell+r\choose r}\G_{{\ell+r\choose r}}=n\})}{\P(\sum j\G_j=n)}
\\
&=&\frac{\P(\G_j=0,\ j\ne{\ell+r\choose r})\P(\sum_\ell{\ell+r\choose r}\G_{{\ell+r\choose r}}=n)}
{\P(\sum j\G_j=n)},
\end{eqnarray*}
where the last identity follows from the fact that the events
$\{\G_j=0,\ j\ne{\ell+r\choose r},\ \ell\ge0\}$ and
$\{\sum_\ell{\ell+r\choose r}\G_{{\ell+r\choose r}}=n\}$
are independent. The same computation yields
\begin{eqnarray*}
&&\Q(\{\la\in A\}\cap\{\la\in P_r(n)\}) \\ &&\qquad =\frac{\P(\G_j=0,\ j\ne{\ell+r\choose r})\P(\{(\G_{{\ell+r\choose r}})\in A\}\cap\{\sum_\ell{\ell+r\choose r}\G_{{\ell+r\choose r}}=n\})}
{\P(\sum j\G_j=n)},
\end{eqnarray*}
and thus
\[
\Pr(\la_{(r)}\in A)=\frac{\P(\{(\G_{{\ell+r\choose r}})\in A\}\cap\{
\sum_\ell{\ell+r\choose r}\G_{{\ell+r\choose r}}=n\})}
{\P(\sum_\ell{\ell+r\choose r}\G_{{\ell+r\choose r}}=n)} \]
Most of the effort is to establish a lower bound on the denominator.   Following  Fristedt \cite{f}
we will show that for a particular choice of the parameter $q=q_n$
there exists a constant $\kappa_r$ such that
\begin{equation}\P(\sum_\ell{\ell+r\choose r}\G_{{\ell+r\choose r}}=n)\ge \kappa_rn^{-\frac{2r+1}{2(r+1)}},\label{lbdd}\end{equation}
for large $n$.
To this end we will choose the value of $q$ which makes
 the expected value of the sum
$$X_n=\sum_{\ell\ge0}{\ell+r\choose r}\G_{{\ell+r\choose r}},$$
asymptotic to $n$ and then we establish a local central limit theorem (CLT) for the normalized random variables $X_n$.
This will show that the probability in question is of order 1 over the standard deviation of $X_n$.
Since $\G$'s are geometric we have
\begin{equation}\E X_n=\sum_{\ell=0}^\infty{\ell+r\choose r}\frac{q^{{\ell+r\choose r}}}{1-q^{{\ell+r\choose r}}},\label{exp}\end{equation}
and since they, in addition, are independent
\begin{equation}\var(X_n)=\sum_{\ell=0}^\infty{\ell+r\choose r}^2\frac{q^{{\ell+r\choose r}}}{\left(1-q^{{\ell+r\choose r}}\right)^2}.
\label{var}\end{equation}
Before proceeding any further, let us note that in subsequent computations, it will be frequently convenient to replace ${\ell+r\choose r}$ by $\ell^r/r!$ in the infinite sums like (\ref{exp}) or (\ref{var}).
This can be done without difficulty because, with our choice of $q$, these sums will grow to infinity  at the faster rate than the individual
terms. Hence,  ignoring the first few terms will not  affect the asymptotic behavior of the sum and for larger $\ell$'s  the approximation ${\ell+r\choose r}\sim\ell^r/r!$ is valid.
 For example, with regard to (\ref{exp}), since the function $uq^u/(1-q^u)$ is decreasing for $u>0$
$$\sum_{\ell=0}^\infty{\ell+r\choose r}\frac{q^{{\ell+r\choose r}}}{1-q^{{\ell+r\choose r}}}\le
\sum_{\ell=1}^\infty\frac{\ell^r}{ r!}\frac{q^{\ell^r/r!}}{(1-q^{\ell^r/r!})}.
$$
On the other hand
\begin{eqnarray*}
&&\sum_{\ell=0}^\infty{\ell+r\choose r}\frac{q^{{\ell+r\choose r}}}{1-q^{{\ell+r\choose r}}}\ge
\sum_{\ell=0}^\infty\frac{(\ell+r)^r}{r!}\frac{q^{(\ell+r)^r/r!}}
{1-q^{(\ell+r)^r/r!}}
\\ && \quad=\sum_{\ell=r}^\infty\frac{\ell^r}{r!}\frac{q^{\ell^r/r!}}
{(1-q^{\ell^r/r!})}\ge\sum_{\ell=1}^\infty\frac{\ell^r}{r!}\frac{q^{\ell^r/r!}}{(1-q^{\ell^r/r!})}-(r-1)\max_{1\le\ell<r}
\left\{\frac{\ell^r}{r!}\frac{q^{\ell^r/r!}}{(1-q^{\ell^r/r!})}\right\}.
\end{eqnarray*}
As we will see, the first term is of order $1/\ln^{(r+1)/r}(1/q)$ and
$\max_{1\le\ell<r}
\left\{\frac{\ell^r}{r!}\frac{q^{\ell^r/r!}}{(1-q^{\ell^r/r!})}\right\}\le
1/\ln(1/q)$. Our choice of $q$ will guarantee that the second term is of lower order than the first one, thus  justifying our claim. Other instances can be treated in a virtually the same manner.
In order to evaluate (\ref{exp}) we consider the sum
$$\sum_{\ell=1}^\infty\frac{\ell^r}{r!}\frac{q^{\ell^r/r!}}{1-q^{\ell^r/r!}}=\sum_{\ell=1}^\infty g(\ell),$$
where
$$g(x)=\frac{x^r}{r!}\frac{q^{x^r/r!}}{1-q^{x^r/r!}}.$$
The function  $g$ is decreasing on a positive half-line.
Therefore,
$$\int_0^\infty g(x)dx\ge\sum_{\ell=1}^\infty g(\ell)\ge \int_1^\infty g(x)dx
\ge \int_0^\infty g(x)dx- g(0+).$$
Hence, an error resulting from replacing the sum by the integral
$$\int_0^\infty g(x)dx=\int_0^\infty \frac{x^r}{r!}\frac{q^{x^r/r!}}{1-q^{x^r/r!}}dx$$
is no more than $g(0+)=1/\ln(1/q)$.
Changing variables to
$q^{x^r/r!}=e^{-u}$ we see that this integral is
$$\frac{(r!)^{1/r}}r\frac1{\ln^{(r+1)/r}(1/q)}\int_0^\infty u^{1/r}\frac{e^{-u}}{1-e^{-u}}du=
\frac{(r!)^{1/r}}r\frac{H_r}{\ln^{(r+1)/r}(1/q)},$$
where $H_r = \zeta(1+r^{-1}) \Gamma(1+r^{-1})$
(see \cite[formula 3.411-7]{gr}).
Thus we can set
$$q=q_n=\exp(-\frac{C}{n^{r/(r+1)}}),$$
where 
$C$ is given by (\ref{Ceq}). 
With that choice, the approximation error between the expectation and $n$ is no more than
$\ln(1/q)=O(n^{r/(r+1)})$.
Asymptotic evaluation of the variance follows the same pattern (with $g(x)$ replaced by $(x^r/r!)^2q^{x^r/r!}/(1-q^{x^r/r!})^2$) and yields
\[
\sigma^2\sim\frac{(r!)^{1/r}}r\frac1{\ln^{(2r+1)/r}(1/q)}\int_0^\infty u^{(r+1)/r}\frac{e^{-u}}{(1-e^{-u})^2}du =\frac{K}{\ln^{\frac{2r+1}r}(1/q)}\sim K n^{\frac{2r+1}{r+1}}
\]
with an error of approximation bounded by $n^{2r/(r+1)}$.
The next step is to establish the local  CLT. Continuing to follow Fristedt, we will first use characteristic functions to establish the CLT.
Let us supress the dependence on $n$ and put
$$X=\sum_{\ell=0}^\infty{\ell+r\choose r}\G_{{\ell+r\choose r}},\quad \mu=\E X,\quad \sigma^2=\var(X),\quad Y=\frac{X-\mu}\sigma.$$
Let
$\phi(t)=\E e^{itY}$
be the characteristic function of $Y$. By independence
\begin{eqnarray*}\log\phi(t)&=&-it\frac\mu\sigma+\sum_{\ell\ge0}\log\left(
\frac{1-q^{{\ell+r\choose r}}}
{1-(q\exp(it/\sigma))^{{\ell+r\choose r}}}
\right)
\\ &=&-it\frac\mu\sigma-\sum_{\ell\ge0}
\log\left(1+
\frac{q^{{\ell+r\choose r}}(1-\exp(i{\ell+r\choose r}t/\sigma))}
{1-q^{{\ell+r\choose r}}}\right).\end{eqnarray*}
Approximating $\log(1+v)$ by $v-v^2/2$ we see that the above series is
asymptotic to
\begin{equation}\sum_{\ell\ge 0}\left\{
\frac{q^{{\ell+r\choose r}}(\exp(i{\ell+r\choose r}t/\sigma)-1)}
{1-q^{{\ell+r\choose r}}}+
\frac{q^{2{\ell+r\choose r}}(1-\exp(i{\ell+r\choose r}t/\sigma))^2}
{2\left(1-q^{{\ell+r\choose r}}\right)^2}
\right\},\label{cltsum}\end{equation}
provided that an  error from approximation is negligible.
But this error is no more than
a constant multiple of
\begin{eqnarray*}&& \sum_{\ell\ge 0}
\frac{q^{3{\ell+r\choose r}}}{\left(1-q^{{\ell+r\choose r}}\right)^3}
\left|1-\exp(i{\ell+r\choose r}t/\sigma)\right|^3
\le
\frac{K|t|^3}{\sigma^3}\sum_{\ell\ge0}{\ell+r\choose r}^3\frac{q^{3{\ell+r\choose r}}}
{\left(1-q^{{\ell+r\choose r}}\right)^3} \\
&&\qquad\sim
\frac{K|t|^3}{\sigma^3}\int_0^\infty\left(\frac{x}{r!}\right)^3
\frac{q^{3x^r/r!}}
{\left(1-q^{x^r/r!}\right)^3}dx=\Theta(n^{-\frac1{2(r+1)}}).
\end{eqnarray*}
Next, we use the approximations
$$1-\exp(i\frac t\sigma{\ell+r\choose r})\sim-
\frac{i{\ell+r\choose r}t}\sigma+
\frac{{\ell+r\choose r}^2t^2}{2\sigma^2}
\quad\textrm{and}\quad
1-\exp(i\frac t\sigma{\ell+r\choose r}
)\sim-\frac{i{\ell+r\choose r}t}\sigma$$
in the first and second expressions in (\ref{cltsum}), respectively. Since the errors are, respectively, of order
$$\frac{{\ell+r\choose r}^3t^3}{\sigma^3}\quad\textrm{and}\quad
\frac{{\ell+r\choose r}t^2}{\sigma^2},$$
 the total errors from approximating these sums are, respectively,
\[
\frac{t^3}{\sigma^3}\sum_{\ell\ge 0}\frac{q^{{\ell+r\choose r}}}{1-q^{{\ell+r\choose r}}}\cdot{\ell+r\choose r}^3
=\Theta\left(\frac1{\sigma^3\ln^{(3r+1)/r}(1/q)}\right)=\Theta(n^{-\frac1{2(r+1)}}),
\]
and
\[K\frac{t^3}{\sigma^3}\int_0^\infty x^{3r}\frac{q^{2x^r/r!}}
{\left(1-q^{x^r/r!}\right)^2}dx
+
K\frac{t^4}{\sigma^4}\int_0^\infty x^{4r}\frac{q^{2x^r/r!}}
{\left(1-q^{x^r/r!}\right)^2}dx
=\Theta(n^{-1/(2(r+1))}),
\]
and are thus negligible.
It follows that
\begin{eqnarray*}
\phi(t)&=&-\frac{it\mu}\sigma+\sum_{\ell\ge0}\frac{q^{{\ell+r\choose r}}}{1-q^{{\ell+r\choose r}}}\cdot\frac{it{\ell+r\choose r}}\sigma
-\sum_{\ell\ge0}\frac{q^{{\ell+r\choose r}}}{1-q^{{\ell+r\choose r}}}\cdot\frac{t^2{\ell+r\choose r}^2}{2\sigma^2}
\\ && \qquad-\sum_{\ell\ge0}\frac{q^{2{\ell+r\choose r}}}{\left(1-q^{{\ell+r\choose r}}\right)^2}\cdot\frac{t^2{\ell+r\choose r}^2}{2\sigma^2}+o(1).
\end{eqnarray*}
The sum of the first two terms is zero while the sum of the remaining two is
$$-\frac{t^2}{2\sigma^2}
\sum_{\ell\ge0}\left\{\frac{q^{{\ell+r\choose r}}}{1-q^{{\ell+r\choose r}}}
\cdot{\ell+r\choose r}^2+\frac{q^{2{\ell+r\choose r}}}{\left(1-q^{{\ell+r\choose r}}\right)^2}\cdot{\ell+r\choose r}^2\right\}.
$$
And again,
\begin{eqnarray*}
& &\sum_{\ell\ge0}\left\{\frac{q^{{\ell+r\choose r}}}{1-q^{{\ell+r\choose r}}}
\cdot{\ell+r\choose r}^2+\frac{q^{2{\ell+r\choose r}}}{\left(1-q^{{\ell+r\choose r}}\right)^2}\cdot{\ell+r\choose r}^2\right\}
\\ & &\qquad\qquad
\sim\frac{(r!)^{1/r}}r\frac1{\ln^{(2r+1)/r}(1/q)}\left\{\int_0^\infty\frac{u^{(r+1)/r}e^{-u}}{1-e^{-u}}du
+\int_0^\infty\frac{u^{(r+1)/r}e^{-2u}}{(1-e^{-u})^2}du\right\}
\\ & & \qquad\qquad =
\frac{(r!)^{1/r}}r\frac1{\ln^{(2r+1)/r}(1/q)}\int_0^\infty\frac{u^{(r+1)/r}e^{-u}(1-e^{-u})+u^{(r+1)/r}e^{-2u}}{(1-e^{-u})^2}du
\\ & & \qquad\qquad =
\frac{(r!)^{1/r}}r\frac1{\ln^{(2r+1)/r}(1/q)}\int_0^\infty\frac{u^{(r+1)/r}e^{-u}}{(1-e^{-u})^2}du,
\end{eqnarray*}
which is the same expression as the one appearing in the computation for $\sigma^2$. Thus, $\forall t\in {\bf R}$
$$\phi(t)\to\exp\left(-\frac{t^2}2\right)$$
which establishes the CLT.

It remains to strenghten it to the local CLT. To this end we will appeal to
\cite[Theorem 2.9]{cs} with
$h_n
=1/\sigma=1/\sigma_n$. We need to find an integrable function $\phi^*$ and a sequence $(\beta_n)$, $\beta_n\to\infty$, such that
\begin{equation}\forall t\qquad\qquad |\phi_n(t)|I(|t|\le\beta_n)\le \phi^*(t),
\label{1st}\end{equation}
and
\begin{equation}\sup_{\beta_n\le|t|\le\pi \sigma_n}|\phi_n(t)|=o(1/\sigma_n).
\label{2nd}\end{equation}
To check (\ref{1st}), supressing a subscript $n$ again, and using the earlier computation  we have
$$\log|\phi(t)|=-\sum_{\ell\ge0}\log\left|\frac{1-q^{{\ell+r\choose r}}\exp(it{\ell+r\choose r}/\sigma)}{1-q^{{\ell+r\choose r}}}\right|$$
and since
$$
\left|\frac{1-q^{{\ell+r\choose r}}\exp(it{\ell+r\choose r}/\sigma)}{1-q^{{\ell+r\choose r}}}\right|
=\left(1+\frac{2q^{{\ell+r\choose r}}(1-\cos(t{\ell+r\choose r}/\sigma))}{\left(1-q^{{\ell+r\choose r}}\right)^2}\right)^{1/2}
$$
we get

\begin{eqnarray*}
\log|\phi(t)|&=&-\frac12\sum_{\ell\ge0}\log
\left(1+\frac{2q^{{\ell+r\choose r}}(1-\cos(t{\ell+r\choose r}/\sigma))}{\left(1-q^{{\ell+r\choose r}}\right)^2}\right) \\ &\le&
-\frac12\sum_{\ell\ge0}\log\left(1+2q^{{\ell+r\choose r}}(1-\cos(t{\ell+r\choose r}/\sigma))\right),
\end{eqnarray*}
by the monotonicity of $\log$. To find a suitable upper bound we increase the right hand side by restricting $\ell$ to the range
\begin{equation}\gamma_r\sigma^{2r/(2r+1)}\le{\ell+r\choose r}\le \sigma^{2r/(2r+1)},\quad0<\gamma_r<1,\label{range}\end{equation}
and replacing the remaining terms by 0's. In that range we have
$$\frac{|t|{\ell+r\choose r}}\sigma\le|t|\sigma^{-1/(2r+1)},$$and
\begin{eqnarray*}
q^{{\ell+r\choose r}}&=&\exp\left(-\frac{C}
{n^{r/(r+1)}}\cdot{\ell+r\choose r}\right)\ge
\exp\left(-\frac{K}
{n^{r/(r+1)}}\cdot\sigma^{\frac{2r}{2r+1}}\right)
\\ &=&\exp(-K)\ge \textrm{const.}
\end{eqnarray*}
Hence, with $\beta_n=\gamma\sigma^{1/(2r+1)}$,
$|t|<\beta_n$ implies $|t|\sigma^{-1/(2r+1)}\le \gamma$, so that
$$2q^{{\ell+r\choose r}}\left(1-\cos\frac t\sigma{\ell+r\choose r}\right)
\ge K
\frac{t^2}{\sigma^2}{\ell+r\choose r}^2,$$
and thus
\begin{eqnarray*}
&&\log\left(1+2q^{{\ell+r\choose r}}\left(1-\cos\frac t\sigma{\ell+r\choose r}\right)\right)\ge K
\frac{t^2}{\sigma^2}{\ell+r\choose r}^2
\\ && \qquad\ge
 K
\frac{t^2}{\sigma^2}\sigma^{\frac{4r}{2r+1}}=Kt^2\sigma^{-\frac2{2r+1}}.\end{eqnarray*}
Since there are at least $K\sigma^{2/(2r+1)}$ $\ell$'s in the range
(\ref{range})
suming within that range gives
$$
\log\left(1+2q^{{\ell+r\choose r}}\left(1-\cos\frac t\sigma{\ell+r\choose r}\right)\right)\ge K\sigma^{\frac2{2r+1}}t^2\sigma^{-\frac2{2r+1}}=Kt^2,$$
and it follows that in order to fulfill (\ref{1st}) we may take
$$\phi^*(t)=\exp(-Kt^2).$$
It remains to verify  (\ref{2nd}).
Pick a $t$ satisfying
$$\gamma\sigma^{1/(2r+1)}\le |t|\le \pi\sigma,$$
and consider the sum
$$
\sum_{\ell\ge0}
\log\left(1+2q^{{\ell+r\choose r}}\left(1-\cos\frac t\sigma{\ell+r\choose r}\right)\right)
$$
Since $\log(1+x)\ge x/5$ for $0\le x\le4$, this sum is at least
$$\frac15
\sum_{\ell\ge0}
q^{{\ell+r\choose r}}\left(1-\cos\frac t\sigma{\ell+r\choose r}\right).
$$
Dealing with the series in the usual way we see that it is asymptotic to
$$\int_0^\infty q^{\frac{x^r}{r!}}\left(1-\cos\left(\frac t\sigma\frac{x^r}{r!}\right)\right)dx,$$
which, using \cite[Formula 3.944-6]{gr}, we find to be asymptotic to
$$\frac{K}{\ln^{1/r}(1/q)}
\left\{1-\frac{\cos\left(\frac1r\arctan(\frac t{\sigma\ln(1/q)})\right)}
{\left(\frac{t^2}{\sigma^2\ln^2(1/q)}+1\right)^{1/(2r)}}\right\}
.
$$
Since $t\ge\gamma\sigma^{1/(2r+1)}$, it follows that
$$\frac{t^2}{\sigma^2\ln^2(1/q)}\ge\gamma^2,$$
and hence, using the relationship between $\sigma$ and $n$ we get
$$\log|\phi(t)|\le -K\frac1{\ln^{1/r}(1/q)}\sim-Kn^{1/(r+1)}\sim-K\sigma^{2/(2r+1)},$$
which implies (\ref{2nd}).
We now appeal to Theorem 2.9 of \cite{cs}
with
$$X_n=\sum_{\ell\ge0}{\ell+r\choose r}\G_{{\ell+r\choose r}},\quad Y_n=\frac{X_n-\mu_n}{\sigma_n},\quad h_n=\frac1{\sigma_n},\quad {\rm and}\quad y_n=\frac{n-\mu_n}{\sigma_n}.$$
Since $n-\mu_n=O(n^{r/(r+1)})$, we have $y_n\to0$
  and  thus
$$\lim_{n\to\infty}n^{\frac{2r+1}{2(r+1)}}\P(X_n=n)
=\lim_{n\to\infty}\frac1{h_n}\P(Y_n=y_n)=\frac1{\sqrt{2\pi}},$$
which clearly implies (\ref{lbdd}).

To complete the proof, let $I_j=\{\G_{{j+r\choose r}}\ge 1\}$. Then, denoting for simplicity a set and its indicator by the same symbol, we have
\begin{eqnarray*}\Pr(|D_{n,r}-\E D_{n,r}|\ge t)
&=&\frac{\P(\{|\sum_j(I_j-{\bf E}I_j)|\ge t\}\cap\{
\sum_\ell{\ell+r\choose r}\G_{{\ell+r\choose r}}=n\})}{
\P(\sum_\ell{\ell+r\choose r}\G_{{\ell+r\choose r}}=n)}
\\&\le& Kn^{\frac{2r+1}{2(r+1)}}\P(|\sum_j(I_j-\E I_j)|\ge t).
\end{eqnarray*}
Since the random variables $(I_j-\E I_j)$ are independent, mean -- zero, and uniformly bounded by 1,  the last probability can be controlled by virtue of Prokhorov's "arcsinh" inequality
(\cite[Theorem 5.2.2(ii)]{s})
$$
\P\left(\sum_j|I_j-\E I_j|\ge t\right)
\le2\exp\left\{-\frac t2{\rm arcsinh}\left(\frac t{2\var(\sum_jI_j)}\right)\right\}
$$
Since
$\var(I_j)=q^{{j+r\choose r}}(1-q^{{j+r\choose r}})$,  by independence we get
\begin{eqnarray*}
\var(\sum_jI_j)&=&\sum_{j\ge0}q^{{j+r\choose r}}(1-q^{{j+r\choose r}})
\sim\frac{(r!)^{1/r}}r\frac1{\ln^{1/r}(1/q)}\int_0^\infty u^{-(r-1)/r}e^{-u}(1-e^{-u})du \\
&=&\frac{(r^!)^{1/r}\left(1-2^{-1/r}\right)\G(1/r)}{r\ln(1/q)}
\sim Kn^{1/(r+1)}.
\end{eqnarray*}
Thus,
selecting $t=t_n$ so that $t_n=o(n^{1/(r+1)})$
 and $t_n^2/n^{1/(r+1)}\to \infty$ not too slow, say,
 $t_n=\Theta(n^{3/(2(r+1))})$, and using the fact
 that ${\rm arcsinh} x\ge Kx$ for $x$ close to 0, we see that
\begin{eqnarray*}
\Pr(|D_{n,r}-\E D_{n,r}|\ge t_n)&\le&
Kn^{(2r+1)/(2(r+1))}\exp\left(-K\frac{t_n^2}{n^{1/(r+1)}}\right)\\ &=&
Kn^{(2r+1)/(2(r+1))}\exp\left(-Kn^{1/(3(r+1))}\right)\to 0.
\end{eqnarray*}
Hence, since  $t_n=o(\E D_{n,r})$, integrating  $D_{n,r,m}/D_{n,r}$ over the set $\{|D_{n,r}-\E D_{n,r}|\le t_n\}$ and its complement yields
$$\E\frac{D_{n,r,m}}{D_{n,r}}= \frac{\E D_{n,r,m}}
{\E D_{n,r}}+O\left(\max\left\{\frac{t_n}{\E D_{n,r}},\P(|D_{n,r}-\E D_{n,r}|>t_n)\right\}\right)=
\frac1{m^{1/r}}+o(1),$$
as desired.

\section{Future work}

One interesting question would be to show how fast the probability
converges. It appears that our proof does come with specific 
rates. The question is therefore to show if they are optimal or not.
Let us note we also showed that the probability that
a random part size of a random partition  into parts 
${{i+r}\choose r}$ of $n$ has multiplicity
at least $m$ is $m^{-1/r}$ when $n\rightarrow\infty$.
Our aim is now to identify the sets S
such  that the probability that
a random part size of a random partition  into parts 
in S
of $n$ has multiplicity
at least $m$ is a constant when $n\rightarrow\infty$.

\end{document}